\documentclass[12pt]{amsart}
\usepackage{amssymb}
\usepackage{graphics}
\usepackage{latexsym}
\usepackage{amsmath}
\usepackage{amssymb,amsthm,amsfonts}
\usepackage{amscd}
\usepackage[arrow, matrix, curve]{xy}
\usepackage{syntonly}
\title[PVHS of CY type and Characteristic Subvariety]{Polarized Variation of Hodge Structures of Calabi-Yau Type and Characteristic Subvarieties Over Bounded
Symmetric Domains}
\author[Mao Sheng]{Mao Sheng}
\author[Kang Zuo]{Kang Zuo}
\address{Department of Mathematics, East China Normal
University, 200062 Shanghhai, P.R. China}
\email{msheng@math.ecnu.edu.cn}
\address{Universit\"{a}t Mainz, Fachbereich 17, Mathematik, 55099 Mainz, Germany}
\email{kzuo@mathematik.uni-mainz.de}
\begin{document}
\theoremstyle{plain}
\newtheorem{thm}{Theorem}[section]
\newtheorem{theorem}[thm]{Theorem}
\newtheorem{lemma}[thm]{Lemma}
\newtheorem{corollary}[thm]{Corollary}
\newtheorem{proposition}[thm]{Proposition}
\newtheorem{addendum}[thm]{Addendum}
\newtheorem{variant}[thm]{Variant}
\theoremstyle{definition}
\newtheorem{construction}[thm]{Construction}
\newtheorem{notations}[thm]{Notations}
\newtheorem{question}[thm]{Question}
\newtheorem{problem}[thm]{Problem}
\newtheorem{remark}[thm]{Remark}
\newtheorem{remarks}[thm]{Remarks}
\newtheorem{definition}[thm]{Definition}
\newtheorem{claim}[thm]{Claim}
\newtheorem{assumption}[thm]{Assumption}
\newtheorem{assumptions}[thm]{Assumptions}
\newtheorem{properties}[thm]{Properties}
\newtheorem{example}[thm]{Example}
\newtheorem{conjecture}[thm]{Conjecture}
\numberwithin{equation}{thm}

\newcommand{\sA}{{\mathcal A}}
\newcommand{\sB}{{\mathcal B}}
\newcommand{\sC}{{\mathcal C}}
\newcommand{\sD}{{\mathcal D}}
\newcommand{\sE}{{\mathcal E}}
\newcommand{\sF}{{\mathcal F}}
\newcommand{\sG}{{\mathcal G}}
\newcommand{\sH}{{\mathcal H}}
\newcommand{\sI}{{\mathcal I}}
\newcommand{\sJ}{{\mathcal J}}
\newcommand{\sK}{{\mathcal K}}
\newcommand{\sL}{{\mathcal L}}
\newcommand{\sM}{{\mathcal M}}
\newcommand{\sN}{{\mathcal N}}
\newcommand{\sO}{{\mathcal O}}
\newcommand{\sP}{{\mathcal P}}
\newcommand{\sQ}{{\mathcal Q}}
\newcommand{\sR}{{\mathcal R}}
\newcommand{\sS}{{\mathcal S}}
\newcommand{\sT}{{\mathcal T}}
\newcommand{\sU}{{\mathcal U}}
\newcommand{\sV}{{\mathcal V}}
\newcommand{\sW}{{\mathcal W}}
\newcommand{\sX}{{\mathcal X}}
\newcommand{\sY}{{\mathcal Y}}
\newcommand{\sZ}{{\mathcal Z}}
\newcommand{\A}{{\mathbb A}}
\newcommand{\B}{{\mathbb B}}
\newcommand{\C}{{\mathbb C}}
\newcommand{\D}{{\mathbb D}}
\newcommand{\E}{{\mathbb E}}
\newcommand{\F}{{\mathbb F}}
\newcommand{\G}{{\mathbb G}}
\newcommand{\HH}{{\mathbb H}}
\newcommand{\I}{{\mathbb I}}
\newcommand{\J}{{\mathbb J}}
\renewcommand{\L}{{\mathbb L}}
\newcommand{\M}{{\mathbb M}}
\newcommand{\N}{{\mathbb N}}
\renewcommand{\P}{{\mathbb P}}
\newcommand{\Q}{{\mathbb Q}}
\newcommand{\R}{{\mathbb R}}
\newcommand{\SSS}{{\mathbb S}}
\newcommand{\T}{{\mathbb T}}
\newcommand{\U}{{\mathbb U}}
\newcommand{\V}{{\mathbb V}}
\newcommand{\W}{{\mathbb W}}
\newcommand{\X}{{\mathbb X}}
\newcommand{\Y}{{\mathbb Y}}
\newcommand{\Z}{{\mathbb Z}}
\newcommand{\id}{{\rm id}}
\newcommand{\rank}{{\rm rank}}
\newcommand{\END}{{\mathbb E}{\rm nd}}
\newcommand{\End}{{\rm End}}
\newcommand{\Hom}{{\rm Hom}}
\newcommand{\Hg}{{\rm Hg}}
\newcommand{\tr}{{\rm tr}}
\newcommand{\Sl}{{\rm Sl}}
\newcommand{\Gl}{{\rm Gl}}
\newcommand{\Cor}{{\rm Cor}}
\newcommand{\proofend}{\hspace*{13cm} $\square$ \\}
\maketitle
\begin{abstract}
In this paper we extend the construction of the canonical polarized
variation of Hodge structures over tube domain considered by B.
Gross in \cite{G} to bounded symmetric domain and introduce a series
of invariants of infinitesimal variation of Hodge structures, which
we call characteristic subvarieties. We prove that the
characteristic subvariety of the canonical polarized variations of
Hodge structures over irreducible bounded symmetric domains are
identified with the characteristic bundles defined by N. Mok in
\cite{M}. We verified the generating property of B. Gross for all
irreducible bounded symmetric domains, which was predicted in
\cite{G}.
\end{abstract}

\section{Introduction}
\footnotetext[1]{This work was supported by the SFB/TR 45 ¡°Periods,
Moduli Spaces and Arithmetic of Algebraic Varieties¡± of the DFG
(German Research Foundation).}

\footnotetext[2]{The first named author is supported by a
Postdoctoral Fellowship in the East China Normal University.}

It has been interesting for long to find a new global Torelli
theorem for Calabi-Yau manifolds, which extends the celebrated
global Torelli theorem for polarized K3 surfaces. The work of B.
Gross \cite{G} is closely related to this problem. In fact, at the
Hodge theoretical level, B. Gross \cite{G} has constructed certain
canonical real polarized variation of Hodge structures (PVHS) over
each irreducible tube domain, and then asked for the possible
algebraic geometrical realizations of them (cf. \cite{G} \S8). In
this paper, we introduce certain invariants, called
\emph{characteristic subvarieties}, which turned out to be
nontrivial obstructions to the realization
problem posed by B. Gross. \\

For a universal family of polarized Calabi-Yau $n$-folds $f: \sX\to
S$, we consider the $\Q$-PVHS $\V$ formed by the primitive middle
rational cohomologies of fibers. Let $(E,\theta)$ be the system of
Hodge bundles associated with $\V$. By the definition of Calabi-Yau
manifold, we have the first property of $(E,\theta)$:
$$
\rank E^{n,0}=1.
$$

The Bogomolov-Todorov-Tian unobstructedness theorem for the moduli
space of Calabi-Yau manifolds gives us the second property of
$(E,\theta)$:
$$
\theta:
T_{S}\stackrel{\simeq}{\longrightarrow}\Hom(E^{n,0},E^{n-1,1}).
$$

On the other hand, we see that the canonical $\R$-PVHSs associated
to tube domains considered by B. Gross in \cite{G} also have the
above two properties (See \cite{G}, Proposition 4.1 and Proposition
5.2). In this paper, we consider the following types of $\C$-PVHS
(See Definition 4.6, \cite{Zuc} for the notion of $\C$-PVHS).
\begin{definition}
Let $\V$ be a $\C$-PVHS of weight $n$ over complex manifold $S$ with
associated system of Hodge bundles $(E,\theta)$. We call $\V$ PVHS
of Calabi-Yau type or {\bf Type I} if $\V$ is a $\R$-PVHS and
$(E,\theta)$ satisfies
\begin{itemize}
    \item $\rank E^{n,0}=1$;
    \item $\theta:
T_{S}\stackrel{\simeq}{\longrightarrow}\Hom(E^{n,0},E^{n-1,1})$.
\end{itemize}
If $\V$ is not defined over $\R$ and $(E,\theta)$ satisfies the
above two properties, then we call $\V$ PVHS of {\bf Type II}.
\end{definition}
For a type I PVHS, one has $\rank E^{0,n}=1$ by the Hodge isometry.
In our case it is not important to emphasize the real structure of a
PVHS. Rather we will simply regard a PVHS of type I or type II
as of Calabi-Yau type.\\

This paper consists of three parts. The first part extends the
construction of PVHS by B. Gross in \cite{G} to each irreducible
bounded symmetric domain. This is a straightforward step. In the
second part, we introduce a series of invariants associated with
the infinitesimal variation of Hodge structures (IVHS) (See the
initiative work of IVHS in \cite{CGGH}). We call these invariants
characteristic subvarieties. Our main result (Theorem
\ref{identification}) identifies the characteristic subvarieties
of the canonical PVHS over an irreducible bounded symmetric domain
with the characteristic bundles defined by N. Mok in \cite{M} (see
\cite{M1} Chapter 6 and Appendix III for a more expository
introduction). The last part of this paper verifies the generating
property predicted by B. Gross (cf. \cite{G} \S5), and describes the
canonical PVHSs over irreducible bounded symmetric domains in
certain detail. \\

In a recent joint work with Ralf Gerkmann \cite{GSZ}, we used the
results of this paper to disprove modularity of the moduli space of
Calabi-Yau 3-folds arising from
eight planes of $\P^3$ in general positions. \\

\textbf{Acknowledgements:} The authors would like to thank Ngaiming
Mok for his explanation of the notion of characteristic bundles, and
thank Eckart Viehweg for his interests and helpful discussions on
this work.

\section{The canonical PVHS over Bounded Symmetric Domain}

Let $D$ be an irreducible bounded symmetric domain, and let $G$ be
the identity component of the automorphism group of $D$. We fix an
origin $0\in D$. Then the isotropy subgroup of $G$ at 0 is a maximal
compact subgroup $K$. By Proposition 1.2.6 \cite{D}, $D$ determines
a special node $v$ of the Dynkin diagram of the simple complex Lie
algebra $\mathfrak{g}^{\C}=Lie(G)\otimes \C$. By the standard theory
on the finite dimensional representations of semi-simple complex Lie
algebras (cf. \cite{FH}), we know that the special node $v$ also
determines a fundamental representation $W$ of $\mathfrak{g}^{\C}$.
By the Weyl's unitary trick, $W$ gives rise to an irreducible
complex representation of $G$. When $D$ is tube domain, the
representation $W$ is exactly the one considered by B. Gross in
\cite{G} and only in this case $W$ admits an $G$-invariant real
form. It is helpful to illustrate the above construction in the
simplest case.
\begin{example}\label{Type A example}
Let $D=SU(p,q)/S(U(p)\times U(q))$ be a type A bounded symmetric
domain. Then
$$
G=SU(p,q),\quad K=S(U(p)\times U(q)),\quad
\mathfrak{g}^{\C}=sl(p+q,\C).
$$
The special node $v$ of the Dynkin diagram of $sl(p+q,\C)$
corresponding to $D$ is the $p$-th node. Let $\C^{p+q}$ be the
standard representation of $Sl(p+q,\C)$. Then the fundamental
representation denoted by $v$ is
$$
W=\bigwedge^{p}\C^{p+q}.
$$
The group $G$ preserves a hermitian symmetric bilinear form $h$ with
signature $(p,q)$ over $\C^{p+q}$. Then $D$ is the parameter space
of the $h$-positive $p$-dimensional vector subspace of $\C^{p+q}$.
By fixing an origin $0\in D$, we obtain an $h$-orthogonal
decomposition
$$
\C^{p+q}=\C^{p}_{+}\oplus \C^{q}_{-}.
$$
The corresponding Higgs bundle to $\W$ is of the form $$
E=\bigoplus_{i+j=n}E^{i,j},$$ where $n=\rm{min}(p,q)$ is the rank of
$D$. The Hodge bundle $E^{n-i,i}$ is the homogeneous vector bundle
determined by, at the origin 0, the irreducible $K$-representation
$$
(E^{n-i,i})_{0}=\bigwedge^{p-i}\C^{p}_{+}\otimes
\bigwedge^{i}\C^{q}_{-}.
$$
\end{example}
Let $\Gamma$ be a torsion free discrete subgroup of $G$, we can
obtain from the representation $W$ the complex local system
$$
\W=W\times_{\Gamma}D
$$
over the locally symmetric variety $X=\Gamma\backslash D$. By the
construction on the last paragraph of \S4, \cite{Zuc}, we know that
$\W$ is a $\C$-PVHS. We denote by $(E,\theta)$ the associated system
of Hodge bundles with $\W$. With similar proofs as those in
Proposition 4.1 and Proposition 5.2 of \cite{G}, or from the
explicit descriptions given in the last section, we have the
following
\begin{theorem}\label{Calabi-Yau like property}
Let $D=G/K$ be an irreducible bounded symmetric domain of rank $n$
and $\Gamma$ be a torsion free discrete subgroup of $G$. Let $\W$ be
the irreducible PVHS over the locally symmetric variety
$X=\Gamma\backslash D$ constructed above. Then $\W$ is a weight $n$
$\C$-PVHS of Calabi-Yau type.
\end{theorem}
Following B. Gross we shall call $\W$ over $X$ in the above theorem
as the canonical PVHS over $X$.

\section{The Characteristic Subvariety and The Main Result}
We start with a system of Hodge bundles
$$
(E=\bigoplus_{p+q=n}E^{p,q},\theta=\bigoplus_{p+q=n}\theta^{p,q})
$$
over a complex manifold $X$ with $\dim E^{n,0}\neq 0$. By the
integrability of Higgs field $\theta$, the $k$-iterated Higgs field
factors as $E\to E\otimes S^k(\Omega_X)$. It induces in turn the
following natural map
$$
\theta^k:S^k(T_X)\to \End(E).
$$
By the Griffiths's horizontal condition, the image of $\theta^k$ is
contained in the subbundle
$$
\bigoplus_{p+q=n}\Hom(E^{p,q},E^{p-k,q+k})\subset \End(E).
$$
We are interested in the projection of $\theta^k$ into the first
component of the above subbundle. Abusing the notation a little bit,
we denote the composition map still by $\theta^k$. That is, we
concern the following map
$$
\theta^k: S^k(T_X)\to \Hom(E^{n,0},E^{n-k,k}).
$$
We have a tautological short exact sequence of analytic coherent
sheaves defined by the iterated Higgs field $\theta^{k}$:
$$
0\to I_k\to S^k(T_X)\stackrel{\theta^k}{\to} J_k\to 0.
$$
We define a sheaf of graded $\sO_X$-algebras $\sJ_{k}$ by putting
$$
\sJ_{k}^{i}=\left\{
  \begin{array}{ll}
    S^{i}\Omega_{X}  & \textrm{if $i< k$}, \\
    S^{i}\Omega_{X}/Im((J_{k})^*\otimes
S^{i-k}\Omega_{X}\stackrel{mult.}{\longrightarrow} S^{i}\Omega_{X})
& \textrm{if $i\geq k$}.
  \end{array}
\right.
$$
\begin{definition}\label{characteristic subvariety}
For $k\geq 0$, we call
$$
C_{k}=Proj(\sJ_{k+1})
$$ the $k$-th
characteristic subvariety of $(E,\theta)$ over $X$.
\end{definition}
By the definition the fiber of a characteristic subvariety over a
point is the zero locus of system of polynomial equations determined
by the Higgs field in the projective tangent space over that point.
In a concrete situation one will be able to calculate some numerical
invariants of the characteristic subvarieties. For example, for a
complete smooth family of hypersurfaces in a projective space, one
can use the Jacobian ring to represent the system of Hodge bundles
associated with the PVHS of the middle dimensional primitive
cohomolgies in a small neighborhood. In \cite{GSZ} one finds such a
calculation in another case.

Because $\theta^{n+1}=0$,
$$
C_{k}=\P(T_X),\ k\geq n,
$$
where $\P(T_X)$ is the projective tangent bundle of $X$. For $0\leq
k \leq n-1$, the natural surjective morphism of graded
$\sO_X$-algebras
$$
\bigoplus_{i=0}^{\infty}S^{i}\Omega_{X}\twoheadrightarrow \sJ_k
$$
gives a proper embedding over $X$,
$$
\xymatrix{
 C_k   \ar[rr]^{\hookrightarrow} \ar[dr]_{p_k}
                &  &   \P(T_X) \ar[dl]^{p}    \\
                & X                }
$$
The next lemma gives a simple criterion to test whether a nonzero
tangent vector at the point $x\in X$ has its image in
$(C_{k})_{x}=p_k^{-1}(x)$.
\begin{lemma}\label{key lemma}
Let $v\in (T_X)_x$ be a non-zero tangent vector at $x$ and $v^{k}\in
(S^k(T_X))_x $ the $k$-th symmetric tensor power of $v$. Then its
image $[v]\in (\P(T_X))_{x}$ lies in $(C_{k-1})_x$ if and only if
$v^k\in (I_k)_{x}$, the stalk of $I_k$ at $x$.
\end{lemma}
\textbf{Proof:} $(C_{k-1})_{x}\subset (\P(T_X))_{x}$ is defined by
the homogeneous elements contained in $((J_k)^*)_x$. Thus $[v]\in
(C_{k-1})_{x}$ if and only for all $f\in ((J_k)^*)_x$, $f([v])=0$.
Now we choose a basis $\{e_1,\cdots,e_m\}$ for $(T_X)_x$ and the
dual
basis $\{e_{1}^{*},\cdots,e_{m}^{*}\}$ for $(\Omega_X)_x$. \\

\textrm{Claim:} $f([v])=0$ if and only $f(v^k)=0$. In the latter, we
consider $f$ as a linear form on
$(S^k(T_X))_x$. \\

\textrm{Proof of Claim:} Let $I=(i_1,\cdots,i_m)$ denote the
multi-index with $i_j\neq 0$ for all $j$, and one puts
$$
I!=i_{1}!\cdots i_{m}!, \quad |I|= i_{1}+\cdots +i_{m}.
$$
We write $v=\sum_{i=1}^{m}a_{i}e_i$ and $f=\sum_{|I|=k}b^{I}e^{I}$.
Then considering $f$ as a polynomial of degree $k$ on $(T_X)_x$, we
have
$$
f(v)=\sum_{|I|=k}b^{I}a^{I}.
$$
On the other hand, we have
$$
v^k=k!\sum_{|I|=k}\frac{1}{I!}a^{I}e^{I},
$$
where $a^{I}=a_{1}^{i_1}\cdots a_{m}^{i_m}$ etc. By
Ex.B.12,\cite{FH} the canonically dual basis of $(S^k(\Omega_X))_x$
to the natural basis $\{e^{I},\ |I|=k\}$ of $(S^k(T_X))_x$ is
$\{\frac{1}{I!}(e^{*})^{I},\ |I|=k\}$. Hence, evaluating $f$ as a
linear form of $(S^k(T_X))_x$ at $v^k$, we obtain
$$
f(v^k)=k!(\sum_{|I|=k}b^{I}a^{I}).
$$
It is clear now that our claim holds. \\

Finally, it is easy to see that $v^k\in (S^k(T_X))_{x}$ lies in
$(I_k)_{x}$ if and only if for all $f\in ((J_k)^*)_x$, considered as
a linear form of $(S^k(T_X))_{x}$, $f(v^k)=0$. Therefore, the lemma
follows.

\proofend

Our main result identifies the characteristic subvarieties of the
canonical PVHS over an irreducible bounded symmetric domain with the
characteristic bundles defined by N. Mok in \cite{M}.
\begin{theorem}\label{identification}
Let $D$ be an irreducible bounded symmetric domain of rank $n$, and
let $(E,\theta)$ be the system of Hodge bundles associated to the
canonical PVHS over $X=\Gamma\backslash D$ as constructed in Theorem
\ref{Calabi-Yau like property}. Then for each $k$ with $1\leq k\leq
n-1$ the $k$-th characteristic subvariety $C_{k}$ of $(E,\theta)$
over $X$ coincides with the $k$-th characteristic bundle $\sS_{k}$
over $X$.
\end{theorem}
By the second property of being of Calabi-Yau type, $C_{0}$ is
always empty. For the self-containedness of this paper, we would
like to describe briefly the notion of characteristic bundles and
refer to Chapter 6 and Appendix III in \cite{M1} for a full
account.\\

The $k$-th characteristic bundle $\sS_k$ over $X=\Gamma\backslash D$
is firstly defined over $D$. It is a projective subvariety of
$\P(T_D)$ and homogeneous under the natural action of automorphism
group $G$ on the projective tangent bundle of $D$. By taking
quotient under the left action of $\Gamma$, one obtains the $k$-th
characteristic bundle over $X$. So it suffices to describe the
construction of characteristic bundle at one point of $D$. At the
origin 0 of $D$, the vectors contained in the fiber $(\sS_{k})_{0}$
are in fact determined by a rank condition. We have the isotropy
representation of $K$ on the tangent space $(T_{D})_{0}$. Fix a
Cartan subalgebra $\mathfrak{h}$ of $\mathfrak{g}$ and choose a
maximal set of strongly orthogonal positive non-compact roots
$$
\Psi=\{\psi_1,\cdots,\psi_{n}\}.
$$
Let $e_{i},1\leq i\leq n,$ be a root vector corresponding to the
root $\psi_{i}$. Then the set $\Psi$ determines a distinguished
polydisk
$$
\triangle^{n}\subset D
$$
passing through the origin 0, and
$$
(T_{\triangle^{n}})_{0}=\sum_{1\leq i \leq n}\C e_{i}\subset
(T_{D})_{0}.
$$
Moreover, for any nonzero element $v\in (T_{D})_{0}$, there exists
an element $k\in K^{\C}$ such that
$$
k(v)=\sum_{1\leq i \leq r(v)}e_{i}.
$$
Such an expression for the vector $v$ is unique and the natural
number $r(v)$ is called the \emph{rank} of $v$. Then, for $1\leq
k\leq n-1$, one defines
$$
(\sS_{k})_{0}=\{[v]\in (\P(T_D))_{0}| 1\leq r(v)\leq k\}.
$$
By the definition, we have a natural inclusion
$$
\sS_{1}\subset \cdots \subset \sS_{n-1}\subset \P(T_D),
$$
We can add two trivially defined characteristic bundles by putting
$$
\sS_0=\emptyset,\quad \sS_n=\P(T_D).
$$
$(\P(T_D))_{0}$ is then decomposed into a disjoint union of
irreducible $K^{\C}$ orbits
$$
(\P(T_D))_{0}=\coprod_{1\leq k\leq n}\{
(\sS_{k})_{0}-(\sS_{k-1})_{0}\}.
$$
\begin{example}
Let $D$ be the type A tube domain of rank $n$. Then
$$
D=SU(n,n)/S(U(n)\times U(n)).
$$
One classically represents $D$ as a space of matrices
$$
D=\{Z\in M_{n,n}(\C)|I_n-\bar{Z}^{t}Z>0\}.
$$
At the origin $0\in D$,
$$
(T_{D})_{0}\simeq M_{n,n}(\C).
$$
The action of
$$
K^{\C}\simeq S(Gl(n,\C)\times Gl(n,\C))
$$
defined by
$$
M\mapsto AMB^{-1}, \ \mathrm{for} \ M\in M_{n,n}(\C) \ \mathrm{and}
\ (A,B)\in Gl(n,\C)\times Gl(n,\C)
$$
gives the isotropy representation of $K^{\C}$ on $(T_{D})_{0}$. Then
the rank of a vector $M\in (T_{D})_{0}$ defined above is just the
rank of $M$ as matrix. Let $(\tilde{\sS}_{k})_{0}$ be the lifting of
$(\sS_k)_0$ in $(T_{D})_{0}$. Therefore, for $1\leq k\leq n-1$, we
have
$$
(\tilde{\sS}_{k})_{0}-(\tilde{\sS}_{k-1})_{0}=S(Gl(n,\C)\times
Gl(n,\C))/P_{k},
$$
where
$$
P_{k}=\{(A,B)\in S(Gl(n,\C)\times Gl(n,\C))| A\left(
                                                 \begin{array}{cc}
                                                   I_{k} & 0 \\
                                                   0 & 0 \\
                                                 \end{array}
                                               \right)=\left(
                                                 \begin{array}{cc}
                                                   I_{k} & 0 \\
                                                   0 & 0 \\
                                                 \end{array}
                                               \right)B
\}.
$$
One can show easily that for $1\leq k\leq n-1$ the dimension of
$(\tilde{\sS}_{k})_{0}-(\tilde{\sS}_{k-1})_{0}$ is $(2n-k)k$.
\end{example}
Now let $D$ be an irreducible bounded symmetric domain of rank $n$,
and let
$$
i: \triangle^n=\triangle_{1}\times \cdots\times
\triangle_{n}\hookrightarrow D
$$
be a polydisc embedding. We are going to study the decomposition of
$i^{*}\W$ into a direct sum of irreducible PVHSs over the polydisc.
The following proposition is a key ingredient in the proof of
Theorem \ref{identification}.
\begin{proposition}\label{decomposition over polydisc}
Let $p_i, 1\leq i \leq n,$ be the projection of the polydisc
$\triangle^n$ into the $i$-th direct factor $\triangle_i$. Then each
irreducible component contained in $i^{*}\W$ is of the form
$$
p_{1}^{*}(\L^{\otimes k_1})\otimes \cdots\otimes
p_{n}^{*}(\L^{\otimes k_n})\otimes \U
$$
with
$$
0\leq k_i\leq 1,\ \textrm{for all}\ \  i,
$$
where $\L$ is the weight 1 PVHS coming from the standard
representation of $Sl(2,\R)$ and $\U$ is a certain unitary factor.
As a consequence, there exists a unique component of the form
$$
p_{1}^{*}\L\otimes \cdots\otimes p_{n}^{*}\L
$$
in $i^{*}\W$ because $\W$ is of Calabi-Yau type.
\end{proposition}
\textbf{Proof:} It is known that the polydisc embedding
$$
i: \triangle^{n}\hookrightarrow D,
$$
determined by a maximal set of strongly orthogonal noncompact roots
$\Psi\subset \mathfrak{h}^{*}$, lifts to a group homomorphism
$$
\phi: Sl(2,\R)^{\times n}\to G.
$$
Our problem is to study the decomposition of $W$ with respect to all
$Sl(2,\R)$ direct factors of $\phi$. \\

We can in fact reduce this to the study of only one direct factor.
This is because a permutation of direct factors can be induced from
an inner automorphism of $G$, which implies the restriction to each
direct factor is isomorphic to each other. Furthermore, we can
assume that the highest root $\tilde{\alpha}$ appears in our chosen
$\Psi$ without loss of generality (cf. \cite{M1} Ch. 5, Proposition 1).\\

Let $s_{\tilde{\alpha}}$ be the distinguished $sl_2$-triple in the
complex simple Lie algebra $\mathfrak{g}^{\C}$ corresponding to
$\tilde{\alpha}$. Let
$$
W=\bigoplus_{\beta\in \Phi} W_{\beta}
$$
be the weight decomposition of $W$ with respect to the Cartan
subalgebra $\mathfrak{h}$. Then by (14.9) \cite{FH}, it is clear
that all irreducible component in $W$ with respect to
$s_{\tilde{\alpha}}$ is contained in
$$
W_{[\beta]}=\bigoplus_{n\in \Z}W_{\beta+n\tilde{\alpha}}.
$$
Let $\mathrm{Conv}(\Phi)$ be the convex hull of $\Phi$, which is a
closed convex polyhedron in $\mathfrak{h}^*$. We put
$$
\partial \Phi=\Phi\cap \mathrm{Conv}(\Phi).
$$
Then for $\beta\in \partial \Phi$, we know by (14.10) \cite{FH} that
the largest component in $W_{[\beta]}$ has dimension equal to
$\beta(H_{\tilde{\alpha}})+1$. Our proof boils
down to showing the following\\

\textrm{Claim:} For all $\beta\in \partial \Phi$, we have
$$
|\beta(H_{\tilde{\alpha}})|\leq 1.
$$

\textrm{Proof of Claim:} We first note that
$$
\displaystyle{ |\beta(H_{\tilde{\alpha}})|=|\frac{2(\beta,
\tilde{\alpha})}{(\tilde{\alpha},\tilde{\alpha})}|}
$$
defines a convex function on $\Phi$. The maximal value will be
achieved for the vertices of $\Phi$, namely, the orbit of highest
weight $\omega$ of $W$ under the Weyl group $\mathrm{W(R)}$. Since
the Weyl group preserves the Killing form, we will show that
$$
|\omega(H_{s(\tilde{\alpha})})|\leq 1, \ \textrm{for all}\ \ s\in
W(R).
$$
The above inequality holds obviously for $s=id$. Let $\alpha_0$ be
the simple root which is the special node determined by $D$ in the
last section. By the definition of special node, in the expression
of $\tilde{\alpha}$ as a linear combination of simple roots the
coefficient before $\alpha_0$ is one (cf. 1.2.5. \cite{D}).
Therefore,
\begin{eqnarray*}
  \omega(H_{\tilde{\alpha}}) &=& \frac{2(\omega,\tilde{\alpha})}{(\tilde{\alpha},\tilde{\alpha})}  \\
   &=&\frac{2(\omega,\alpha_0)}{(\tilde{\alpha},\tilde{\alpha})}\\
   &=&\frac{(\alpha_0,\alpha_0)}{(\tilde{\alpha},\tilde{\alpha})}\\
   &=&1.
\end{eqnarray*}
Now we have to separate the exceptional cases from the ongoing proof
because of the complicated description of the Weyl group in the
exceptional cases. In the following, we use the same notation as the
appendix of \cite{B}. Let $\{\varepsilon_1,\cdots,\varepsilon_{l}\}$
be the standard basis of the Euclidean space $\R^{l}$, and $\sigma$
denotes a permutation of index.\\

Type $A_{l-1}$: The highest root
$\tilde{\alpha}=\varepsilon_1-\varepsilon_{l}$. The Weyl group
permutes the basis elements. All fundamental weights
$$\omega_{i}=\sum_{j=1}^{i}\varepsilon_j-\frac{i}{l+1}\sum_{j=1}^{l}\varepsilon_j, 1\leq
i\leq l-1$$ correspond to a special node. Then
\begin{eqnarray*}
  |\omega_{i}(H_{s(\tilde{\alpha})})| &=& |(\omega_{i},s(\tilde{\alpha}))|  \\
  &=&|(\sum_{j=1}^{i}\varepsilon_j,\varepsilon_{\sigma(1)}-\varepsilon_{\sigma(l)})|\\
   &\leq &1.
\end{eqnarray*}

Type $B_l$: The highest root
$\tilde{\alpha}=\varepsilon_1+\varepsilon_{2}$. The Weyl group
permutes the basis elements, or acts by $\varepsilon_{i}\mapsto \pm
\varepsilon_{i}$. The first fundamental weight
$\omega_1=\varepsilon_1$ corresponds to the special node. Then
\begin{eqnarray*}
  |\omega_{1}(H_{s(\tilde{\alpha})})| &=& |(\omega_{1},s(\tilde{\alpha}))|  \\
   &=&|(\varepsilon_1,\pm \varepsilon_{\sigma(1)}\pm \varepsilon_{\sigma(2)})|\\
   &\leq &1.
\end{eqnarray*}

Type $C_l$: The highest root $\tilde{\alpha}=2\varepsilon_1$. The
Weyl group permutes the basis elements, or acts by
$\varepsilon_{i}\mapsto \pm \varepsilon_{i}$. The last fundamental
weight $\omega_l=\sum_{i=1}^{l}\varepsilon_i$ corresponds to the
special node.Then
\begin{eqnarray*}
  |\omega_{l}(H_{s(\tilde{\alpha})})| &=& |\frac{1}{2}(\omega_{l},s(\tilde{\alpha}))|  \\
   &=&|\frac{1}{2}(\sum_{i=1}^{l}\varepsilon_i,\pm 2\varepsilon_{\sigma(1)})|\\
   &=&1.
\end{eqnarray*}

Type $D_{l}$: The highest root
$\tilde{\alpha}=\varepsilon_1+\varepsilon_2$. The Weyl group
permutes the basis elements, or acts by $\varepsilon_{i}\mapsto (\pm
1)_{i} \varepsilon_{i}$ with $\prod_{i}(\pm 1)_{i}=1$. We have three
special nodes in this case. It suffices to check
$\omega_{1}=\varepsilon_1$ and
$\omega_{l}=\frac{1}{2}(\sum_{i=1}^{l}\varepsilon_i)$. For
$\omega_{1}$, we have
\begin{eqnarray*}
  |\omega_{1}(H_{s(\tilde{\alpha})})| &=& |(\omega_{1},s(\tilde{\alpha}))|  \\
   &=&|(\varepsilon_1,\pm \varepsilon_{\sigma(1)}\pm \varepsilon_{\sigma(2)})|\\
   &\leq&1.
\end{eqnarray*}
For $\omega_{l}$, we have
\begin{eqnarray*}
  |\omega_{l}(H_{s(\tilde{\alpha})})| &=& |(\omega_{l},s(\tilde{\alpha}))|  \\
   &=&|\frac{1}{2}(\sum_{i=1}^{l}\varepsilon_i,\pm \varepsilon_{\sigma(1)}\pm \varepsilon_{\sigma(2)})|\\
   &\leq &1.
\end{eqnarray*}
Now we treat with the exceptional cases. In the following, we shall
compute the largest value of $|\beta(H_{\tilde{\alpha}})|$ among all
weights $\beta$ in $\Phi$. The results will particularly imply the claim.\\

Type $E_6$: Let $\{\alpha_1,\cdots, \alpha_6\}$ be the set of simple
roots of simple Lie algebra of type $E_6$ and
$\{\omega_1,\cdots,\omega_6\}$ be the fundamental weights. The
highest root is then
$$
\tilde{\alpha}=\alpha_1+2\alpha_2+2\alpha_3+3\alpha_4+2\alpha_5+\alpha_6.
$$
A 6-tuple $(a_1,\cdots a_6)$ denotes the weight
$\beta=\sum_{i=1}^{6}a_{i}\omega_{i}$. There are two special nodes
in this case and it suffices to study either of them. The following
table lists all elements of $\Phi$ for the fundamental
representation $\omega_1$:
$$
\begin{array}{cccc}
  (1,0,0,0,0,0)&(-1,0,1,0,0,0)&(0,0,-1,1,0,0)&(0,1,0,-1,1,0) \\
  (0,1,0,0,-1,1)&(0,-1,0,0,1,0)&(0,1,0,0,0,-1)&(0,-1,0,1,-1,1)\\
  (0,0,1,-1,0,1)&(0,-1,0,1,0,-1)&(1,0,-1,0,0,1)&(0,0,1,-1,1,-1)\\
  (1,0,-1,0,1,-1)&(0,0,1,0,-1,0)&(-1,0,0,0,0,1)&(1,0,-1,1,-1,0)\\
  (-1,0,0,0,1,-1)&(1,1,0,-1,0,0)&(-1,0,0,1,-1,0)&(1,-1,0,0,0,0)\\
  (-1,1,1,-1,0,0)&(0,1,-1,0,0,0)&(-1,-1,1,0,0,0)&(0,-1,-1,1,0,0)\\
  (0,0,0,-1,1,0)&(0,0,0,0,-1,1)&(0,0,0,0,0,-1).&
\end{array}
$$
For an element $(a_1,\cdots a_6)$ in the above table, we have
\begin{eqnarray*}
  \beta(H_{\tilde{\alpha}}) &=& (\beta,\tilde{\alpha})  \\
   &=&(\sum_{i=1}^{6}a_{i}\omega_{i},\alpha_1+2\alpha_2+2\alpha_3+3\alpha_4+2\alpha_5+\alpha_6 )\\
   &= & a_{1}+2a_{2}+2a_{3}+3a_{4}+2a_{5}+a_{6}.
\end{eqnarray*}
According to this formula, it is straightforward to compute that the
largest value of $|\beta(H_{\tilde{\alpha}})|$ is equal to one.\\

Type $E_7$: Let $\{\alpha_1,\cdots,\alpha_7\}$ be the set of simple
roots of simple Lie algebra of type $E_7$ and
$\{\omega_1,\cdots,\omega_7\}$ be the fundamental weights. We can
choose the maximal set $\Psi$ of the strongly orthogonal noncompact
roots to be
$$
\{\psi_1=2\alpha_1+2\alpha_2+3\alpha_3+4\alpha_4+3\alpha_5+2\alpha_6+\alpha_7,
\psi_2=\alpha_2+\alpha_3+2\alpha_4+2\alpha_5+2\alpha_6+\alpha_7,
\psi_3=\alpha_7\}.
$$
It is simpler to use $\psi_3$ to verify our statement, instead of
$\psi_1$ which is the highest root $\tilde{\alpha}$. As in the last
case, a 7-tuple $(a_1,\cdots a_7)$ denotes the weight
$\beta=\sum_{i=1}^{7}a_{i}\omega_{i}$. The following table lists all
elements of $\Phi$ for the fundamental representation $\omega_7$:
$$
\begin{array}{cccc}
  (0,0,0,0,0,0,1)&(0,0,0,0,0,1,-1)&(0,0,0,0,1,-1,0)&(0,0,0,1,-1,0,0)\\
  (0,1,1,-1,0,0,0)&(1,1,-1,0,0,0,0)&(0,-1,1,0,0,0,0)&(1,-1,-1,1,0,0,0)\\
  (-1,1,0,0,0,0,0)&(1,0,0,-1,1,0,0)&(-1,-1,0,1,0,0,0)&(1,0,0,0,-1,1,0) \\
  (-1,0,1,-1,1,0,0)&(1,0,0,0,0,-1,1)&(0,0,-1,0,1,0,0)&(-1,0,1,0,-1,1,0) \\
  (1,0,0,0,0,0,-1)&(0,0,-1,1,-1,1,0)&(-1,0,1,0,0,-1,1)&(0,1,0,-1,0,1,0) \\
  (0,0,-1,1,0,-1,1)&(-1,0,1,0,0,0,-1)&(0,1,0,-1,1,-1,1)&(0,0,-1,1,0,0,-1)\\
  (0,-1,0,0,0,1,0)&(0,1,0,0,-1,0,1)&(0,1,0,-1,1,0,-1)&(0,-1,0,0,1,-1,1)\\
  (0,1,0,0,-1,1,-1)&(0,-1,0,1,-1,0,1)&(0,-1,0,0,1,0,-1)&(0,1,0,0,0,-1,0) \\
  (0,0,1,-1,0,0,1)&(0,-1,0,1,-1,1,-1)&(1,0,-1,0,0,0,1)&(0,0,1,-1,0,1,-1) \\
 (0,-1,0,1,0,-1,0)&(1,0,-1,0,0,1,-1)&(0,0,1,-1,1,-1,0)&(-1,0,0,0,0,0,1) \\
 (1,0,-1,0,1,-1,0)&(0,0,1,0,-1,0,0)&(-1,0,0,0,0,1,-1)&(1,0,-1,1,-1,0,0) \\
  (-1,0,0,0,1,-1,0)&(1,1,0,-1,0,0,0)&(-1,0,0,1,-1,0,0)&(1,-1,0,0,0,0,0)\\
  (-1,1,1,-1,0,0,0)&(0,1,-1,0,0,0,0)&(-1,-1,1,0,0,0,0)&(0,-1,-1,1,0,0,0)\\
 (0,0,0,-1,1,0,0)&(0,0,0,0,-1,1,0)&(0,0,0,0,0,-1,1)&(0,0,0,0,0,0,-1).
\end{array}
$$
Then for an element $(a_1,\cdots a_7)$ in the above table, we have
\begin{eqnarray*}
  \beta(H_{\alpha_7}) &=& (\beta, \alpha_7)  \\
   &=&(\sum_{i=1}^{7}a_{i}\omega_{i},\alpha_7 )\\
   &= & a_{7}.
\end{eqnarray*}
It is straightforward to see the largest value of $|
\beta(H_{\alpha_7})|$ is one. This completes the whole proof.

\proofend

We can now proceed to prove our main result.

\textbf{Proof of Theorem \ref{identification}:} It suffices to prove
the isomorphism over $D$, and we obtain the claimed isomorphism by
taking quotient under the left action of $\Gamma$. Since the
constructions on both sides are $G$-equivariant, it is enough to
show the isomorphism at the origin of $D$. Over the origin 0, we
have the adjoint action of $K^{\C}$ on the holomorphic tangent space
$(T_{D})_{0}$ and the dual action on $(\Omega_{D})_{0}$. Since the
Higgs field of a locally homogeneous VHS is $G$-equivariant, for
each $k$, $(J_{k})_0\subset (S^{k}\Omega_{D})_{0}$ is
$K^{\C}$-invariant. This implies $C_k$ is $K^{\C}$-invariant. So we
can obtain a decomposition of $(\P(T_D))_{0}$ into disjoint union of
$K^{\C}$ orbits as follows:
$$
(\P(T_D))_{0}=\coprod_{1\leq k\leq n} \{(C_{k})_{0}-(C_{k-1})_{0}\}.
$$
For $1\leq k \leq n$, we put
$$
v_k=e_1+\cdots+e_k.
$$
It is clear that
$$
v_k\in (\sS_{k})_{0}-(\sS_{k-1})_{0}
$$
and $K^{\C}(v_{k})=(\sS_{k})_{0}-(\sS_{k-1})_{0}$. Next we make the
following\\

\textrm{Claim:} For $1\leq k\leq n$,
$$
v_k\in (C_{k})_{0}-(C_{k-1})_{0}.
$$
This claim implies the inclusion of $K^{\C}$ orbit for each $k$
$$
(\sS_{k})_{0}-(\sS_{k-1})_{0}\subset (C_{k})_{0}-(C_{k-1})_{0},
$$
and hence the equality for each $k$. Therefore, for $1\leq k\leq
n-1$,
$$
(C_{k})_{0}=(\sS_{k})_{0}.
$$

\textrm{Proof of Claim:} Let $\gamma_k: \triangle \to D$ be the
composition map
$$
\triangle \stackrel{\rm{diag.}}{\longrightarrow
}\triangle_1\times\cdots\times\triangle_k\hookrightarrow
\triangle^n\stackrel{i}{\hookrightarrow }D.
$$
Obviously, for a suitable basis element $u$ of $(T_{\triangle})_0$
one has $(d\gamma_{k})_{0}(u)=v_k$. By Proposition
\ref{decomposition over polydisc}, we have the decomposition of PVHS
$$
\gamma_{k}^{*}\W=\L^{\otimes k}\otimes \U\oplus \V^{'}
$$
where $\U$ is a unitary factor and $\V^{'}$ is a PVHS with width
$\leq k-1$. Let $(E,\theta)$ be the system of Hodge bundles
corresponding to $\W$, and for $v\in (T_{D})_0$, we denote by
$$
\theta_{v}: E_{0}\to E_{0}
$$
the action of the Higgs field $\theta$ on the bundle $E$ at the
origin 0 along the tangent direction $v$. Since
$$
\theta_{v_{k}}(E)=\theta_{u}(\gamma_{k}^{*}E),
$$
we see that
$$
(\theta_{v_k})^k\neq 0,\quad (\theta_{v_k})^{k+1}=0.
$$
Together with Lemma \ref{key lemma}, one easily sees that the claim
holds.

\proofend

\section{Enumeration of Canonical PVHS over Irreducible Bounded Symmetric Domain and the Generating Property of Gross}
Let $(E,\theta)$ be a system of Hodge bundles over $X$. We use the
same notation as that in previous sections. We note that
$$
I=\bigoplus_{k\geq 1}I_{k}
$$
forms a graded ideal of the symmetric algebra
$$
Sym(T_X)=\bigoplus_{k\geq 0}S^{k}T_{X}.
$$
It is trivial to see that
$$
I_{k}=S^{k}(T_{X}),\ \textrm{for}\ k\geq n+1.
$$
In \cite{G} \S5, Gross suspected if $I$ is generated by $I_2$ for
the canonical PVHS over an irreducible tube domain. We can assert
this generating property for the canonical PVHS over an irreducible
bounded symmetric domain in general.
\begin{theorem}\label{generating property}
We use the same notation as Theorem \ref{Calabi-Yau like property}.
Then the graded ideal $I$, formed by the kernel of iterated Higgs
field, is generated by the degree 2 graded piece $I_2$. That is, the
multiplication map
$$
I_{2}\otimes S^{k-2}(T_X)\to I_{k}
$$
is surjective for all $k\geq 2$.
\end{theorem}
It suffices to prove the surjectivity for $k\leq n+1$ where
$n=\rank(D)$. In fact, for $k\geq n+2$, we have
\begin{eqnarray*}
  I_{2}\otimes S^{n-1}(T_{X})\otimes (T_{X})^{\otimes
k-n-1} & \twoheadrightarrow& I_{n+1}\otimes (T_{X})^{\otimes k-n-1} \\
    &=&  S^{n+1}(T_{X})\otimes (T_{X})^{\otimes k-n-1} \\
    &\twoheadrightarrow& S^{k}(T_{X})=I_{k}.
\end{eqnarray*}
By the integrality of Higgs field, the above surjective map factors
through $I_{2}\otimes S^{k-2}(T_{X})$. As the proof of Theorem
\ref{identification}, we can work on the level of bounded symmetric
domain and prove the statement at the origin as $K$-representations.
The theorem will be proved case by case. In the classical case, we
shall also describe the system of Hodge bundles $(E,\theta)$
associated with the Calabi-Yau like PVHS $\W$ using the Grassmannian
description of classical symmetric domain. \\

Let $D$ be an irreducible bounded symmetric domain. By fixing an
origin of $D$, we obtain an equivalence of categories of homogeneous
vector bundles and finite dimensional complex representations of
$K$. Since $K$ has one dimensional center, a finite dimensional
complex $K$-representation is written as $\C(l)\otimes V$ where $V$
is a representation of the semisimple part $K'$ of $K$ and is
determined by the induced action of the complexified Lie algebra
$\mathfrak{k'^{\C}}$. In the following, the same notation for a
$K$-representation and the corresponding homogeneous vector bundle
will be used when the context causes no confusion. All the
isomorphisms are isomorphisms between homogeneous bundles. A highest
weight representation of $sl(n,\C)$ will be denoted interchangeably
by $\Gamma_{a_1,\cdots,a_{n-1}}$ and $\SSS_{\lambda}(\C^{n})$ (cf.
\cite{FH} \S15.3).

\subsection{Type $A$}
The irreducible bounded symmetric domain of type A is
$D^{I}_{p,q}=G/K$ where
$$
G=SU(p,q),\quad K=S(U(p)\times U(q)).
$$
Let $V=\C^{p+q}$ be a complex vector space equipped with a Hermitian
symmetric bilinear form $h$ of signature $(p,q)$. Then $D^{I}_{p,q}$
parameterizes the dimension $p$ complex vector subspaces $U\subset
V$ such that
$$
h|_{U}: U\times U\to \C
$$
is positive definite. This forms the tautological subbundle
$S\subset V\times D$ of rank $p$ and denote by $Q$ the tautological
quotient bundle of rank $q$. We have the natural isomorphism of
holomorphic vector bundles
\begin{eqnarray}\label{equation1}
  T_{D^{I}_{p,q}} &\simeq & \Hom(S,Q).
\end{eqnarray}
The standard representation $V$ of $G$ gives rise to a weight 1 PVHS
$\V$ over $D^{I}_{p,q}$, and its associated Higgs bundle
$$
F=F^{1,0}\oplus F^{0,1},\quad\eta=\eta^{1,0}\oplus \eta^{0,1}
$$
is determined by
$$
F^{1,0}=S,\quad F^{0,1}=Q,\quad \eta^{0,1}=0,
$$
and $\eta^{1,0}$ is defined by the above isomorphism. The canonical
PVHS is
$$
\W=\bigwedge^{p}\V
$$
and its associated system of Hodge bundles $(E,\theta)$ is then
$$
(E,\theta)=\bigwedge^{p}(F,\eta).
$$
Since
$$
\mathfrak{k'^{\C}}=sl(p,\C)\oplus sl(q,\C),
$$
by Schur's lemma, a finite dimensional irreducible complex
representation of $\mathfrak{k'^{\C}}$ is of the form
$$
\Gamma_{a_1,\cdots,a_{p-1}}\otimes \Gamma'_{b_1,\cdots,b_{q-1}}.
$$
We put $V_1=\C^{p}$ to be the representation space
$\Gamma_{0,\cdots,0,1}$ of $sl(p,\C)$ and $V_2=\C^{q}$ the
representation space $\Gamma'_{0,\cdots,0,1}$ of $sl(q,\C)$. In the
remaining subsection, we shall assume $p\leq q$ in order to simplify
the notations in the argument.
\begin{lemma}\label{formula A}
We have isomorphism
$$
T_{D^{I}_{p,q}}\simeq V_1\otimes V_2.
$$
Then, for $k\geq 2$, we have isomorphism
$$
S^{k}(T_{D^{I}_{p,q}})\simeq
\bigoplus_{\lambda}\SSS_{\lambda}(V_1)\otimes \SSS_{\lambda}(V_2),
$$
where $\lambda$ runs through all partitions of $k$ with at most $p$
rows. Under this isomorphism, the $k$-th iterated Higgs field for
$k\leq p$,
$$
\theta^{k}: S^{k}(T_{D^{I}_{p,q}})\to  \Hom(E^{p,0},E^{p-k,k})
$$
is identified with the projection map onto the irreducible component
$$
\bigoplus_{\lambda}\SSS_{\lambda}(V_1)\otimes
\SSS_{\lambda}(V_2)\twoheadrightarrow \SSS_{\lambda^{0}}(V_1)\otimes
\SSS_{\lambda^{0}}(V_2),
$$
where $\lambda^{0}=(1,\cdots,1)$.
\end{lemma}
\textbf{Proof:} By the isomorphism \ref{equation1}, we have
isomorphism
$$
T_{D^{I}_{p,q}}\simeq V_1\otimes V_2.
$$
The formula in Ex.6.11 \cite{FH} gives the decomposition of
$S^{k}(V_1\otimes V_2)$ with respect to $sl(p,\C)\oplus sl(q,\C)$:
$$
S^{k}(V_1\otimes V_2)=\bigoplus_{\lambda}\SSS_{\lambda}(V_1)\otimes
\SSS_{\lambda}(V_2),
$$
where $\lambda$ runs through all partitions of $k$ with at most $p$
rows. Since the center of $K$ acts on $(T_{D^{I}_{p,q}})_{0}$
trivially, it acts on $(S^{k}(T_{D^{I}_{p,q}}))_{0}$ trivially too.
Hence the second isomorphism of the statement follows. For the last
statement, it suffices to show $\theta^k$ is a non-zero map because
$\Hom(E^{p,0},E^{p-k,k})$ is irreducible. But this follows directly
from the definition of the Higgs field $\theta$ as $p$-th wedge
power of $\eta$. The lemma is proved.

\proofend

From the lemma, we know that
$$
\theta^2 \simeq pr: \Gamma_{0,\cdots,2} \otimes
\Gamma'_{0,\cdots,2}\oplus \Gamma_{0,\cdots,0,1,0} \otimes
\Gamma'_{0,\cdots,0,1,0}\to \Gamma_{0,\cdots,0,1,0} \otimes
\Gamma'_{0,\cdots,0,1,0}.
$$
So by definition,
$$
I_{2}\simeq  \Gamma_{0,\cdots,0,2}\otimes \Gamma'_{0,\cdots,0,2}.
$$

{\bf Proof of Theorem \ref{generating property} for Type A:} Now we
proceed to prove that $I_{2}\otimes S^{k-2}(T_{D^{I}_{p,q}})$
generates $I_{k}$. By the above lemma and the Formula 6.8 \cite{FH},
we have
\begin{eqnarray*}
 I_{2}\otimes S^{k-2}(T_{D^{I}_{p,q}}) &\simeq& S^2(V_1)\otimes S^2(V_2)\otimes S^{k-2}(V_1\otimes V_2) \\
    &\simeq & \bigoplus_{\mu}(S^2(V_1)\otimes \SSS_{\mu}(V_1))\otimes (S^2(V_2)\otimes \SSS_{\mu}(V_2))  \\
    &=& \bigoplus_{\mu}[(\bigoplus_{\nu_{\mu}^{1}}\SSS_{\nu_{\mu}^{1}}(V_1))\otimes (\bigoplus_{\nu_{\mu}^{2}}\SSS_{\nu_{\mu}^{2}}(V_2)) ]\\
    &=& \bigoplus_{\mu}\bigoplus_{\nu_{\mu}^{1},\nu_{\mu}^{2}}(\SSS_{\nu_{\mu}^{1}}(V_1)\otimes \SSS_{\nu_{\mu}^{2}}(V_2)), \\
\end{eqnarray*}
where $\mu$ runs through all partitions of $k-2$ with at most $p$
rows, and for a fixed $\mu$, $\nu_{\mu}^{i},i=1,2$ runs through
those Young diagrams by adding two boxes to different columns of the
Young diagram of $\mu$. Let $\lambda$ be a Young diagram
corresponding to a direct factor of $I_k$ under the isomorphism in
the above lemma. Since
$$
\bigoplus_{\mu}\bigoplus_{\nu_{\mu}}(\SSS_{\nu_{\mu}}(V_1)\otimes
\SSS_{\nu_{\mu}}(V_2)) \subset
\bigoplus_{\mu}\bigoplus_{\nu_{\mu}^{1},\nu_{\mu}^{2}}(\SSS_{\nu_{\mu}^{1}}(V_1)\otimes
\SSS_{\nu_{\mu}^{2}}(V_2)),
$$
it is enough to show that $\lambda$ can be obtained by a Young
diagram $\mu$ by adding two boxes to different columns of $\mu$.
Actually, by the above lemma the partition $\lambda$ of $I_{k}$ has
the property that, either for some $1\leq i_{0}\leq p-1$,
$$
\lambda_{i_{0}}>\lambda_{i_{0}+1}\geq 1,
$$
or for some $1\leq i_{0}\leq p$,
$$
2 \leq \lambda_1=\cdots=\lambda_{i_{0}}>\lambda_{i_{0}+1}\geq 0.
$$
In the first case, we can choose $\mu$ as
$$
\mu_{i}=\left\{
  \begin{array}{ll}
    \lambda_{i}-1  & \textrm{if $i=i_0, i_{0}+1$}, \\
     \lambda_{i}
& \textrm{otherwise}.
  \end{array}
\right.
$$
In the second case, we choose $\mu$ as
$$
\mu_{i}=\left\{
  \begin{array}{ll}
    \lambda_{i}-2  & \textrm{if $i=i_0$}, \\
     \lambda_{i}
& \textrm{otherwise}.
  \end{array}
\right.
$$
The proof of Theorem \ref{generating property} in the type A case is
therefore completed.

\proofend

\subsection{Type $B$, Type $D^{\R}$} For $n\geq 3$, we let
$$
G=Spin(2,n),\quad K=Spin(2)\times_{\mu_2} Spin(n).
$$
Then $D^{IV}_{n}=G/K$ is the bounded symmetric domain of type B when
$n$ is odd, of type $D^{\R}$ when $n$ even. Let $(V_{\R},Q)$ be a
real vector space of dimension $n+2$ equipped with a symmetric
bilinear form of signature $(2,n)$. Then $D^{IV}_{n}$ is one of
connected components parameterizing all $Q$-positive two dimensional
subspace of $V_{\R}$. In order to see clearer the complex structure
of $D^{IV}_{n}$, we complexify $(V_{\R},Q)$, to obtain
$(V=V_{\R}\otimes \C,Q)$. Then it is known that $D^{IV}_{n}$ is an
open submanifold of the quadratic hypersurface defined by $Q=0$ in
$\P(V)\simeq \P^{n+1}$, which is just the compact dual of
$D^{IV}_{n}$. For a $Q$-isotropic line $L\subset V$, we define its
polarization hyperplane to be
$$
P(L)=\{v\in V| Q(L,v)=0\}.
$$
So for each point of $D^{IV}_{n}$, we obtain a natural filtration of
$V$ by
$$
L\subset P(L)\subset V.
$$
Varying the points on $D^{IV}_{n}$, the above filtration yields a
filtration of homogeneous bundles
$$
S\subset P(S)\subset V\times D^{IV}_{n}.
$$
On the other hand, we have a commutative diagram
$$
\begin{CD}
T_{D^{IV}_{n}} @>\simeq >>  \Hom(L,\frac{P(L)}{L}) \\
@V\cap VV @VV\cap V  \\
T_{\P(V),[L]} @>\simeq >> \Hom(L,\frac{V}{L}),
\end{CD}
$$
whose top horizontal line gives the isomorphism of tangent bundle
\begin{eqnarray}\label{equation2}
  T_{D^{IV}_{n}}  &\simeq & \Hom(S,\frac{P(S)}{S}).
\end{eqnarray}
We also notice that $Q$ descends to a non-degenerate bilinear form
on $\frac{P(L)}{L}$, so that we have a natural isomorphism
\begin{eqnarray}\label{equation2'}
 \big(\frac{P(S)}{S}\big)^*  &\simeq & \frac{P(S)}{S}.
\end{eqnarray}
Now we put
$$
E^{2,0}=S,\quad E^{1,1}=\frac{P(S)}{S},\quad E^{0,2}=\frac{V\times
D^{IV}_{n}}{P(S)},
$$
and
\begin{eqnarray*}
  \theta^{2,0}: E^{2,0} &\to & E^{1,1}\otimes \Omega_{D^{IV}_{n}}, \\
  \theta^{1,1}: E^{1,1} &\to & E^{0,2}\otimes \Omega_{D^{IV}_{n}}
\end{eqnarray*}
are determined by the isomorphisms \ref{equation2} and
\ref{equation2'}, and $\theta^{0,2}=0$. The Higgs bundle
$$
(E=\bigoplus_{p+q=2}E^{p,q},\theta=\bigoplus_{p+q=2}\theta^{p,q})
$$
is the associated system of Hodge bundles with the canonical PVHS
$\W$. \\

Let $m=[\frac{n}{2}]$ be the rank of $\mathfrak{so(n)}$, and
$\Gamma_{a_1,\cdots,a_{m}}$ denotes a highest weight representation
of $\mathfrak{so(n)}$. In terms of this notation, we have
$$
E^{2,0}\simeq\C(-2)\otimes \Gamma_{0,\cdots,0},\quad
E^{1,1}\simeq\C\otimes \Gamma_{1,0,\cdots,0},\quad
E^{0,2}\simeq\C(2)\otimes \Gamma_{0,\cdots,0}.
$$
The following easy lemma makes Theorem \ref{generating property} in
the cases of type B and type $D^{\R}$ clear.
\begin{lemma}\label{formula B}
We have isomorphisms
\begin{eqnarray*}
  T_{D^{IV}_{n}} &\simeq& \C(2)\otimes \Gamma_{1,0,\cdots,0},  \\
  S^{2}(T_{D^{IV}_{n}}) &\simeq& \C(4)\otimes \Gamma_ {2,0,\cdots,0}\oplus \C(4)\otimes \Gamma_{0,\cdots,0}, \\
  I_{2} &\simeq& \C(4)\otimes \Gamma_ {2,0,\cdots,0}, \\
  I_{2}\otimes T_{D^{IV}_{n}}&\simeq& \C(6)\otimes\Gamma_ {3,0,\cdots,0} \oplus \C(6)\otimes \Gamma_ {1,0,\cdots,0}\oplus \C(6)\otimes \Gamma_ {1,1,0,\cdots,0},   \\
  S^{3}(T_{D^{IV}_{n}}) &\simeq&\C(6)\otimes \Gamma_ {3,0,\cdots,0}\oplus \C(6)\otimes \Gamma_
  {1,0,\cdots,0}.
\end{eqnarray*}
\end{lemma}

\subsection{Type $C$} We fix $n\geq 2$. Let
$$
G=Sp(2n,\R),\quad K=U(n).
$$
Then $D^{III}_{n}=G/K$ is the bounded symmetric domain of type C.
$D^{III}_{n}$ is known as the Siegel space of degree $n$. Let
$(V_{\R},\omega)$ be a real vector space of dimension $2n$ equipped
with a skew symmetric bilinear form $\omega$. As before, we denote
also by $(V,\omega)$ the complexification, and
$h(u,v)=i\omega(u,\bar{v}) $ defines a hermitian symmetric bilinear
form over $V$. Then $D^{III}_{n}$ parameterizes the maximal
$\omega$-isotropic and $h$-positive complex subspaces of $V$. The
standard representation $V$ of $G$ gives a weight 1 $\R$-PVHS $\V$
over $D^{III}_{n}$. Let $(F,\eta)$ be the associated Higgs bundle
with $\V$. Then $F^{1,0}$ is simply the tautological subbundle over
$D^{III}_{n}$ and $F^{0,1}$ is the $h$-orthogonal complement of
$F^{1,0}$. Clearly, we have a natural embedding of bounded symmetric
domains
$$
\iota: D^{III}_{n}\hookrightarrow D^{I}_{n,n}.
$$
It induces a commutative diagram:
$$
\begin{CD}
T_{D^{III}_{n}} @>\simeq >> S^{2}(F^{0,1}) \\
@V\cap VV @VV\cap V  \\
\iota^{*}(T_{D^{I}_{n,n}}) @>\simeq >> (F^{0,1})^{\otimes 2}.
\end{CD}
$$
The Higgs field $\eta^{1,0}$ is defined by the composition of maps
\begin{eqnarray}\label{equation3}
T_{D^{III}_{n}}  &\simeq&  S^{2}(F^{0,1}) \hookrightarrow
(F^{0,1})^{\otimes 2}\simeq \Hom(F^{1,0},F^{0,1}).
\end{eqnarray}
The canonical PVHS $\W$ is the unique weight $n$ sub-PVHS of
$\bigwedge^{n}(\V)$. In fact, we have a decomposition of $\R$-PVHS
$$
\bigwedge^{n}(\V)=\W\oplus \V',
$$
where $\V'$ is a weight $n-2$ $\R$-PVHS. Therefore the corresponding
Higgs bundle $(E,\theta)$ to $\W$ is a sub-Higgs bundle of
$\bigwedge^{n}(F,\eta)$.\\

Let $V_1=(F^{0,1})_{0}$ be the standard representation of $K$. It is
straightforward to obtain the following
\begin{lemma}\label{formula C}
We have isomorphism
$$
T_{D^{III}_{n}}\simeq \SSS_{(2)}(V_1).
$$
Then, for $k\geq 2$, we have isomorphism
$$
S^{k}(T_{D^{III}_{n}})\simeq
\bigoplus_{\lambda}\SSS_{\lambda}(V_{1}),
$$
where $\lambda=\{\lambda_1,\cdots,\lambda_{l}\}$ runs through all
partitions of $2k$ with each $\lambda_{i}$ even and $l\leq n$. Under
this isomorphism, for $k\leq n$, the $k$-th iterated Higgs field
$\theta^{k}$ is identified with the projection map onto the
irreducible component $\SSS_{\lambda^{0}}(V_{1})$ where
$\lambda^{0}=(2,\cdots,2)$.
\end{lemma}

{\bf Proof of Theorem \ref{generating property} for Type C:} By the
last lemma, we know that
$$
\theta^2\simeq pr: \SSS_{(4)}(V_{1})\oplus \SSS_{(2,2)}(V_{1})\to
\SSS_{(2,2)}(V_{1})
$$
and then $I_2\simeq \SSS_{(4)}(V_{1})$. Applying the Formula
6.8,\cite{FH} to decompose $I_2\otimes S^{k-2}(T_{D^{III}_{n}})$, we
obtain
\begin{eqnarray*}
 I_2\otimes S^{k-2}(T_{D^{III}_{n}})&\simeq& \SSS_{(4)}(V_{1}) \otimes \bigoplus_{\mu}\SSS_{\mu}(V_{1}) \\
    &\simeq & \bigoplus_{\mu}(\SSS_{(4)}(V_{1}) \otimes \SSS_{\mu}(V_{1}))  \\
    &=& \bigoplus_{\mu}[(\bigoplus_{\nu_{\mu}}\SSS_{\nu}(V_1))],
\end{eqnarray*}
where $\mu$ runs through all partitions of $2(k-2)$ with the
property as that in Lemma \ref{formula C}, and for a fixed $\mu$,
$\nu_{\mu}$ runs through those Young diagrams by adding four boxes
to different columns of the Young diagram $\mu$. The partition
$\lambda$ of an irreducible component in $I_k$ is of the form
$$
\lambda_1\geq \cdots \geq
\lambda_{s}>\lambda_{s+1}=\cdots=\lambda_{l}\geq 2.
$$
We may then take $\mu$ to be
$$
\mu_{i}=\left\{
  \begin{array}{ll}
    \lambda_{i}-2  & \textrm{if $i=s, l$}, \\
     \lambda_{i}
& \textrm{otherwise}.
  \end{array}
\right.
$$
Then we define $\nu_{\mu}$ by adding two boxes to $\lambda_{s}$ and
$\lambda_{l}$ in $\mu$ respectively to obtain the starting
$\lambda$. Therefore Theorem \ref{generating property} in the type C
case is proved.

\proofend

\subsection{Type $D^{\HH}$}
For $n\geq 3$, we let
$$
G=SO^{*}(2n),\quad K=U(n).
$$
Then $D^{II}_{n}=G/K$ is the bounded symmetric domain of type
$D^{\HH}$. We recall that
$$
G\simeq \{M\in Sl(2n,\C)|MI_{n,n}M^{*}=
                          I_{n,n}, MS_{n}M^{\tau}=
                          S_{n}
\},
$$
where $I_{n,n}$ denotes the matrix
$\left(
\begin{array}{cc}
I_n & 0 \\
                              0 & -I_n \\
                            \end{array}
                          \right)$
and $S_n$ denotes the matrix
$\left(
                            \begin{array}{cc}
                               0& I_n \\
                              I_n & 0 \\
                            \end{array}
                          \right)$.
Let $(V,h,S)$ be a complex vector space of dimension $2n$ equipped
with a hermitian symmetric form $h$ and symmetric bilinear form $S$,
where, under the identification $V\simeq \C^{2n}$, $h$ is defined by
the matrix $I_{n,n}$ and $S$ is defined by the matrix $S_n$. Then
$D^{II}_{n}$ parameterizes all $n$-dimensional $S$-isotropic and
$h$-positive complex subspaces of $V$. The standard representation
$V$ of $G$ determines a weight 1 PVHS $\V$. Its associated Higgs
bundle $(F,\eta)$ is determined in a similar manner as type C case.
Namely, $F^{1,0}$ is simply the tautological subbundle and $F^{0,1}$
is its $h$-orthogonal complement. The natural embedding
$$
\iota': D^{II}_{n}\hookrightarrow D^{I}_{n,n}
$$
induces a commutative diagram:
$$
\begin{CD}
T_{D^{II}_{n}} @>\simeq >> \bigwedge^{2}(F^{0,1}) \\
@V\cap VV @VV\cap V  \\
\iota'^{*}(T_{D^{I}_{n,n}}) @>\simeq >> (F^{0,1})^{\otimes 2},
\end{CD}
$$
and the Higgs field $\eta^{1,0}$ is induced by the composition of
maps
\begin{eqnarray}\label{equation4}
T_{D^{II}_{n}}&\simeq&  \bigwedge^{2}(F^{0,1}) \hookrightarrow
(F^{0,1})^{\otimes 2}\simeq \Hom(F^{1,0},F^{0,1}).
\end{eqnarray}
The canonical PVHS $\W$ comes from a half spin representation. We
write the corresponding Higgs bundle as
$$
(E=\bigoplus_{p+q=[\frac{n}{2}]}E^{p,q},\theta=\bigoplus_{p+q=[\frac{n}{2}]}\theta^{p,q}).
$$
Then the Hodge bundle is
$$
E^{p,q}= \bigwedge^{n-2q}F^{1,0},
$$
and the Higgs field $\theta^{p,q}$ is induced by the natural wedge
product map
$$
\bigwedge^{2}F^{0,1}\otimes \bigwedge^{2q}F^{0,1}\to
\bigwedge^{2q+2}F^{0,1}.
$$
While type $D^{\HH}$ case enjoys many similarity with type C case,
there is one difference we would like to point out. That is, the
canonical PVHS $\W$ is not a sub-PVHS of $\bigwedge^{n}\V$. In fact,
the PVHS $\bigwedge^{n}\V$ is the direct sum of two irreducible
PVHSs. One of them, say $\V'$, has
$$\bigwedge^{n}(F^{1,0})\otimes \bigwedge^{0}(F^{0,1})\simeq
(\bigwedge^{n}(F^{1,0}))^{\otimes 2}$$ as the first Hodge bundle.
For this irreducible $\V'$, we have an inclusion of PVHS
$$
\V'\subset Sym^{2}(\W).
$$
Let $V_1=(F^{0,1})_{0}$ be the dual of the standard representation
of $K$. It is straightforward to obtain the following
\begin{lemma}\label{formula D}
We have isomorphism
$$
T_{D^{II}_{n}} \simeq \SSS_{(1,1)}(V_1).
$$
Then, for $k\geq 2$, we have isomorphism
$$
S^{k}(T_{D^{II}_{n}})\simeq
\bigoplus_{\lambda}\SSS_{\lambda}(V_{1}),
$$
where $\lambda=\{\lambda_1,\cdots,\lambda_{l}\}$ runs through all
partitions of $2k$ with $l\leq n$ and each entry of the conjugate
$\lambda'$ of $\lambda$ even. Under this isomorphism, for $k\leq
[\frac{n}{2}]$, the $k$-th iterated Higgs field $\theta^{k}$ is
identified with the projection map onto the irreducible component
$\SSS_{\lambda^{0}}(V_{1})$ with $\lambda^{0}=(k,k)$.
\end{lemma}
By the above lemma, we know that
$$
\theta^2\simeq pr: \SSS_{(1,1,1,1)}(V_{1})\oplus
\SSS_{(2,2)}(V_{1})\to \SSS_{(2,2)}(V_{1}).
$$
Thus we have isomorphism $I_2\simeq \SSS_{(1,1,1,1)}(V_{1})$. By
Formula 6.9,\cite{FH}, we have
\begin{eqnarray*}
 I_2\otimes S^{k-2}(T_{D^{II}_{n}}) &
 \simeq & \SSS_{(1,1,1,1)}(V_{1}) \otimes \bigoplus_{\mu}\SSS_{\mu}(V_{1}) \\
    &\simeq & \bigoplus_{\mu}(\SSS_{(1,1,1,1)}(V_{1}) \otimes \SSS_{\mu}(V_{1}))  \\
    &=& \bigoplus_{\mu}[(\bigoplus_{\nu_{\mu}}\SSS_{\nu}(V_1))],
\end{eqnarray*}
where $\mu$ runs through all partitions of $2(k-2)$ with the
property as that in Lemma \ref{formula D}, and for a fixed $\mu$,
$\nu_{\mu}$ runs through those Young diagrams by adding four boxes
to different rows of the Young diagram of $\mu$. We observe that we
are in the conjugate case of that of type C. Theorem \ref{generating
property} in type $D^{\HH}$ case follows easily.

\subsection{Type $E$}
There are two exceptional irreducible bounded symmetric domains. We
first discuss the $E_6$ case. In this case,
$$
G=E_{6,2},\quad K=U(1)\times_{\mu_4}Spin(10).
$$
Then $D^{V}=G/K$ is a 16-dimensional bounded symmetric domain of
rank 2. There are two special nodes in the Dynkin diagram of $E_6$.
But they induce isomorphic bounded symmetric domains. We take the
first node so that the fundamental representation corresponding to
this special node is $W_{27}$. Let
$(E=\oplus_{p+q=2}E^{p,q},\theta)$ be the corresponding Higgs bundle
to $\W$. Then we have isomorphism
$$
E^{2,0}\simeq \C(-2),\quad E^{1,1}\simeq \C \otimes
\Gamma_{0,0,0,1,0},\quad E^{0,2}\simeq \C(2)\otimes
\Gamma_{1,0,0,0,0}.
$$
Furthermore, it is straightforward to obtain the following
\begin{lemma}\label{formula E6}
We have following isomorphisms:
\begin{eqnarray*}
  T_{X} &\simeq & \C(2)\otimes \Gamma_{0,0,0,1,0},  \\
  S^{2}(T_{X}) &\simeq & \C(4)\otimes \Gamma_{0,0,0,2,0}\oplus \C(4)\otimes \Gamma_{1,0,0,0,0}, \\
  I_{2} &\simeq& \C(4)\otimes \Gamma_{0,0,0,2,0}, \\
  I_{2}\otimes T_{X}&\simeq& \C(6)\otimes \Gamma_{0,0,0,3,0}\oplus \C(6)\otimes \Gamma_{1,0,0,1,0}\oplus \C(6)\otimes \Gamma_{0,0,1,1,0},  \\
  S^{3}(T_{X}) &\simeq &\C(6)\otimes \Gamma_{0,0,0,3,0}\oplus \C(6)\otimes
  \Gamma_{1,0,0,1,0}.
\end{eqnarray*}
\end{lemma}
We continue to discuss the remaining case, which has already
appeared in \cite{G}. Let
$$
G=E_{7,3},\quad  K=U(1)\times_{\mu_3}E_{6}.
$$
Then $D^{VI}=G/K$ is of dimension 27 and rank 3. We refer the reader
to \S4 \cite{G} for the description of Hodge bundles. The lemma
corresponding to Lemma \ref{formula E6} is the following
\begin{lemma}\label{formula E7}
We have the following isomorphisms:
\begin{eqnarray*}
  T_{X} &\simeq & \C(2)\otimes \Gamma_{1,0,0,0,0,0},  \\
  S^{2}(T_{X}) &\simeq& \C(4)\otimes \Gamma_{2,0,0,0,0,0} \oplus \C(4)\otimes \Gamma_{0,0,0,0,0,1}, \\
  I_{2} &\simeq& \C(4)\otimes \Gamma_{2,0,0,0,0,0}, \\
  I_{2}\otimes T_{X}&\simeq& \C(6)\otimes \Gamma_{1,0,0,0,0,1}\oplus \C(6)\otimes \Gamma_{3,0,0,0,0,0}\oplus \C(6)\otimes \Gamma_{1,0,1,0,0,0},  \\
  S^{3}(T_{X}) &\simeq& \C(6)\otimes \Gamma_{1,0,0,0,0,1}\oplus \C(6)\otimes
  \Gamma_{3,0,0,0,0,0}\oplus \C(6)\otimes \Gamma_{0,0,0,0,0,0},\\
  I_{3} &\simeq& \C(6)\otimes \Gamma_{1,0,0,0,0,1}\oplus \C(6)\otimes
  \Gamma_{3,0,0,0,0,0}, \\
  I_{2}\otimes S^{2}(T_{X}) &\simeq&\C(8)\otimes \Gamma_{4,0,0,0,0,0}\oplus \C(8)\otimes \Gamma_{2,0,0,0,0,1}\oplus \C(8)\otimes
  \Gamma_{0,0,0,0,0,2}\\
  &\phantom{\simeq}& \oplus \C(8)\otimes \Gamma_{2,0,1,0,0,0}\oplus \C(8)\otimes
   \Gamma_{0,0,2,0,0,0} \oplus \C(8)\otimes \Gamma_{0,0,1,0,0,1}\\
 &\phantom{\simeq}& \oplus \C(8)\otimes \Gamma_{1,0,0,0,0,0}\oplus \C(8)\otimes
   \Gamma_{1,1,0,0,0,0} \oplus \C(8)\otimes \Gamma_{2,0,0,0,0,1},\\
   S^{4}(T_{X}) &\simeq& \C(8)\otimes \Gamma_{4,0,0,0,0,0}\oplus \C(8)\otimes \Gamma_{2,0,0,0,0,1}\oplus \C(8)\otimes
  \Gamma_{0,0,0,0,0,2} \\
   &\phantom{\simeq}& \oplus \C(8)\otimes \Gamma_{1,0,0,0,0,0}.
\end{eqnarray*}
\end{lemma}
Lemma \ref{formula E6} and Lemma \ref{formula E7} make it clear that
the generating property of Gross also holds for the exceptional
cases. Then the proof of Theorem \ref{generating property} is
completed.

\end{document}